\documentclass{amsart}
\usepackage{geometry} 
\usepackage{url}
\usepackage{amsthm}
\usepackage{amssymb}
\usepackage{latexsym}
\usepackage{amsfonts}
\usepackage{amsmath}
\usepackage{amstext}
\usepackage{accents}
\usepackage{cancel}
\usepackage{eucal}
\usepackage{amscd}
\usepackage{hyperref}
\usepackage[shortalphabetic]{amsrefs}
\usepackage{tikz}

\begin{document}

\title[High order flow of curves]{High order curvature flows of plane curves with generalised Neumann boundary conditions} 
\author[J. McCoy]{James McCoy*} \thanks{* Corresponding author}
\address{Priority Research Centre Computer Assisted Research Mathematics and Applications, School of Mathematical and Physical Sciences, University of Newcastle, Australia and Okinawa Institute for Science and Technology Graduate University, Japan}
\email{James.McCoy@newcastle.edu.au}
\author[G. Wheeler]{Glen Wheeler}
\address{Institute for Mathematics and its Applications, University of Wollongong}
\email{glenw@uow.edu.au}
\author[Y. Wu]{Yuhan Wu}
\address{Institute for Mathematics and its Applications, University of Wollongong}
\email{yw120@uowmail.edu.au}
 \thanks{The research of the first author was supported by Discovery Project DP180100431 of the Australian Research Council.  Part of this work was completed while the first author was a Visiting Professor at the Okinawa Institute for Science and Technology.  The research of the third author was supported by a University of Wollongong Faculty of Engineering and Information Sciences Postgraduate research scholarship.  The authors are grateful for the support provided by these facilities.}
\begin{abstract}
We consider the parabolic polyharmonic diffusion and $L^2$-gradient flows of the $m$-th arclength derivative of curvature for regular closed curves evolving with generalised Neumann boundary conditions.  In the polyharmonic case, we prove that if the curvature of the initial curve is small in $L^2$, then the evolving curve converges exponentially in the $C^\infty$ topology to a straight horizontal line segment.  The same behaviour is shown for the $L^2$-gradient flow provided the energy of the initial curve is sufficiently small.  In each case the smallness conditions depend only on $m$.
\end{abstract}

\maketitle

\section{Introduction} \label{S:intro}
\newtheorem{main1}{Theorem}[section]
\newtheorem{main2}[main1]{Theorem}

Higher order geometric evolution problems have interesting practical applications that have motivated increasing attention in recent years to their theoretical behavior.  
As fourth order examples for evolving curves we have the curve diffusion flow and the $L^2$-gradient flow of the
elastic energy, and for surfaces the corresponding surface diffusion and Willmore flows.
Flows of higher even order than four have been less thoroughly investigated,
but they and their elliptic counterparts are well-motivated given, for example, applications in
computer design, where higher order equations  allow more
flexibility in terms of prescribing boundary conditions \cite{LX}.  Such
equations have also found applications in medical imaging \cite{UW}.  In \cite{MPW}, the first and second author together with Parkins considered the sixth order geometric triharmonic flow for closed surfaces while Parkins and the third order considered in \cite{PW16} even order flows of closed, planar curves.  There the flows were of general even order of polyharmonic form as we will consider here, but we will also in this article consider the $L^2$-gradient flows of the $m$-th arclength derivative of curvature for general $m\in \mathbb{N}\cup \left\{ 0\right\}$.  Our work here generalises \cite{MWW} where we considered the $L^2$-gradient flow for the energy
$$\int_\gamma k_s^2 ds \mbox{;}$$
$k_s$ denotes the first derivative of curvature with respect to the arc
length parameter $s$.  Our work is also the arbitrary even order generalisation of \cite{WW} by the third author and V-M Wheeler, where the fourth order curve diffusion and elastic flow of curves between parallel lines were investigated. Other relevant works on fourth order flow of curves with boundary
conditions are \cites{DLP14, DP14, L12}.  By way of comparison, closed curves without boundary evolving by higher order equations  have been more thoroughly studied; see for example \cites{DKS02, EGBMWW14,
GI99, PW16, W13}.

Let $\gamma_0:\left[ -1, 1\right] \rightarrow \mathbb{R}^2$ be a (suitably) smooth embedded or immersed regular curve whose ends $\gamma_0\left( \pm 1\right)$ meet orthogonally two parallel lines $\eta_{\pm}$ separated by distance $d_0$.  In this article we are interested in one-parameter families of curves $\gamma\left( \cdot, t\right)$ satisfying either the polyharmonic curve diffusion flow
\begin{equation} \label{E:PHCD}
  \frac{\partial \gamma}{\partial t} = \left( -1\right)^{m+1} k_{s^{2m+2}} \nu
\end{equation}
or the flow
\begin{equation} \label{E:theflow}
  \frac{\partial \gamma}{\partial t} = \left[ \left( -1\right)^{m+1} k_{s^{2m+2}} - \sum_{j=1}^m \left( -1\right)^j k\, k_{s^{m+j}} k_{s^{m-j}} - \frac{1}{2} k\, k_{s^m}^2\right] \nu \mbox{,}
\end{equation}
with $\gamma\left( \cdot, 0 \right) := \gamma_0$ and generalised Neumann boundary conditions.  The equation \eqref{E:theflow} corresponds to the $L^2$-gradient flow for the energy
\begin{equation} \label{E:E}
  \int_\gamma k_{s^m}^2 ds \mbox{.}
\end{equation}
Above $k_{s^m}$ denotes the $m$th iterated derivative of curvature with respect to the arc
length parameter $s$; $\nu$ is the smooth choice of unit normal such that the above flows are parabolic in the generalised sense.  As discussed in \cite{Pthesis} for example, \eqref{E:PHCD} can also be considered as a gradient flow in an appropriate corresponding Sobolev space.  The `generalised Neumann boundary conditions' we assume are not the most general possible for either flow, but they are mathematically a natural choice: we take classical Neumann boundary conditions, as shown in Figure 1, together with no curvature flux on the boundary ($k_s\left( \pm 1, t\right)=0$) and we additionally specify that all odd derivatives of curvature up to order $2m+1$ are equal to zero on the boundary.  For each of our flows, induction arguments analogous to those in \cite[Lemma 2.6]{WW} then show that all higher odd curvature derivatives are also equal to zero on the boundary under the flow, so we have for every $\ell \in \mathbb{N}$, as long a solution to the flow equation exists,
\begin{equation} \label{E:bdy}
k_{s^{2\ell-1}}\left( \pm1, t\right) = 0 \mbox{.}
\end{equation}

\noindent \emph{Remark:} If we think of the corresponding higher order elliptic ordinary differential equation, an order $2m+4$ equation should normally have $2m+4$ boundary conditions for a unique solution.  These correspond to the classical Neumann condition and all odd curvature derivatives up to order $2m+1$ equal to zero at $x=\pm 1$.  On the other hand, from the point of view of partial differential equations, it is more natural to think of each pair of boundary requirements at $\pm 1$ as one condition, giving a total of $m+2$ conditions.  We will use the latter description of the number of boundary conditions throughout the article. 

\begin{figure}
\begin{center}
\begin{tikzpicture}

\draw (0,0) node[below] {$\eta_-$} --(0,5) ;
\draw (4,0) node[below] {$\eta_+$}--(4,5);

\draw [ cyan] (0,2.2) to (0.2, 2.2) to [out=0,in=180] (1.5,2.8) to [out=0,in=180] (2.3,2.2) node[below] {$\gamma$} to [out=0,in=-180] (3.8, 3) to (4,3);

\draw (0,2.4)-- (0.2,2.4) -- (0.2,2.2);

\draw (4,2.8)-- (3.8,2.8) -- (3.8,3);

\draw[line width=1pt, -latex, blue] (4,3) --node[auto]{$e$} (5,3);
\draw[line width=1pt, -latex, blue] (4,3) -- node[auto]{$\nu$}(4,4);

\draw [blue] (0,1.4) to  [out=0,in=180] (1.8,1.4) node[] {$ \ \ \ d_0$};

\draw [blue] (0,1.4) -- (0.2, 1.5);
\draw [blue] (0,1.4) -- (0.2, 1.3);

\draw [blue] (2.2,1.4) to [out=0,in=-180] (4,1.4);

\draw [blue] (4, 1.4)--(3.8, 1.3);

\draw [blue] (4, 1.4)--(3.8, 1.5);
\end{tikzpicture}
\ \

Figure 1
\end{center}
\end{figure}

Throughout this article we use $\omega$ to denote the \emph{winding number}, defined here by
\[
\omega := \frac{1}{2\pi} \int_{\gamma} k \, ds \mbox{.}
\]

A simple calculation shows that under quite general flows with Neumann boundary conditions on parallel lines, the winding number is constant \cite[Lemma 2.5]{WW}.\\[8pt]

Our main results are as follows:

\begin{main1} \label{T:main1}
Let $\gamma_0: \left[ -1, 1\right] \rightarrow \mathbb{R}^2$ be a smooth embedded or immersed curve with $\omega=0$, whose ends meet the parallel lines $\eta_{\pm}$ with $m+2$ generalised Neumann boundary conditions as described above.  If the curvature $\kappa$ of $\gamma_0$ is sufficiently small in $L^2$, that is
$$\int \kappa^2 ds \leq \varepsilon$$
for some $\varepsilon>0$ depending only on $m$, then there exists a smooth solution $\gamma :\left[ 0, \infty\right) \rightarrow \mathbb{R}^2$ to \eqref{E:PHCD} with $\gamma\left( \cdot, 0 \right) = \gamma_0$.  The solution $\gamma$ is unique up to parametrisation, smooth and converges exponentially to a horizontal line segment whose distance from $\gamma_0$ is finite.\\
\end{main1}

\begin{main2} \label{T:main2}
Let $\gamma_0: \left[ -1, 1\right] \rightarrow \mathbb{R}^2$ be a smooth embedded or immersed curve with $\omega=0$, whose ends meet the parallel lines $\eta_{\pm}$ with $m+2$ generalised Neumann boundary conditions as described above.  If $\gamma_0$ has sufficiently small energy, that is
$$\int \kappa_{s^m}^2 ds \leq \varepsilon$$
for some $\varepsilon>0$ depending only on $m$, then there exists a smooth solution $\gamma :\left[ 0, \infty\right) \rightarrow \mathbb{R}^2$ to the $L^2$-gradient flow for \eqref{E:E} with $\gamma\left( \cdot, 0 \right) = \gamma_0$.  The solution $\gamma$ is unique up to parametrisation, smooth and converges exponentially to a horizontal line segment whose distance from $\gamma_0$ is finite.\\
\end{main2}

We remark that local existence of a smooth regular solution $\gamma:\left[ -1, 1\right] \times \left( 0, T\right)\rightarrow \mathbb{R}^2$ to each of the above problems for some $T>0$ is standard.  Such solutions are unique up to parametrisation.  If $\gamma_0$ satisfies appropriate compatibility conditions, then the solution is smooth on $\left[ 0, T\right)$.  It is possible to consider such flows with less smooth initial data but we will not do so in this article.  An overview of the procedure for proving short-time existence in this setting is given in \cite[Theorem 2.1]{MWW}.  In particular, the initial curve $\gamma_0$ need not be a graph over the horizontal line segment although of course, the later results show that the solution does indeed eventually become so.\\  

The structure of the rest of this article is as follows.  In Section \ref{S:prelim} we will state fundamental analytical tools that will be used in the analysis of each of our flow problems.  We also give the general structure of some evolution equations that is useful in both cases.  In Section \ref{S:PHCD} we prove Theorem \ref{T:main1}, the polyharmonic curve diffusion case.  The proof has the same structure as the proof of Theorem \ref{T:main2}, thus illustrating the key ideas, however the estimates are much simpler to establish.  In Section \ref{S:setup} we take the normal variation of the energy \eqref{E:E} to obtain the corresponding $L^2$-gradient
flow \eqref{E:theflow}, then we prove Theorem \ref{T:main2}.  
 
\section{Preliminaries} \label{S:prelim}
\newtheorem{PSW}{Lemma}[section]
\newtheorem{interp}[PSW]{Proposition}
\newtheorem{evlneqns}[PSW]{Lemma}

We begin with the following standard result for functions of one variable.
\begin{PSW}[Poincar\'{e}-Sobolev-Wirtinger (PSW) inequalities] \label{T:PSW}
 Suppose $f:\left[ 0, L \right] \rightarrow \mathbb{R}$, $L>0$ is absolutely continuous.  
 \begin{itemize}
   \item If $\int_0^L f \,ds =0$ then
 $$\int_0^L f^2 ds \leq \frac{L^2}{\pi^2} \int_0^L f_s^2 ds \mbox{ and } \left\| f\right\|_\infty^2 \leq \frac{2L}{\pi} \int_0^L f_s^2 ds \mbox{.}$$
  \item Alternatively, if $f\left( 0\right) = f\left( L \right) =0$ then
 $$\int_0^L f^2 ds \leq \frac{L^2}{\pi^2} \int_0^L f_s^2 ds \mbox{ and } \left\| f\right\|_\infty^2 \leq \frac{L}{\pi} \int_0^L f_s^2 ds \mbox{.}$$
 \end{itemize}
 \end{PSW}
 
 To state the next interpolation inequality we will use, we first need to set up
 some notation.  For normal tensor fields $S$ and $T$ we denote by $S \star T$
 any linear combination of $S$ and $T$.  In our setting, $S$ and $T$ will be
 simply curvature $k$ or its arclength derivatives.  Denote by $P_n^m\left(
 k\right)$ any linear combination of terms of type $\partial_s^{i_1} k \star
 \partial_s^{i_2}k \star \ldots \star \partial_s^{i_n}k$ where $m=i_1 + \ldots+
 i_n$ is the total number of derivatives.
 
 It is convenient to use the following scale-invariant norms: we define
 $$\left\| k \right\|_{\ell, p} := \sum_{i=0}^\ell \left\| \partial_s^i k \right\|_p$$
 where
 $$\left\| \partial_s^i k \right\|_p = L^{i+1-\frac{1}{p}} \left( \int \left| \partial_s^i k \right|^p ds\right)^{\frac{1}{p}} \mbox{.}$$
 
 The following interpolation inequality for closed curves appears in
 \cite{DKS02}; for our setting with boundary we refer to \cite{DP14}.
 
 \begin{interp} \label{T:int}
 Let $\gamma: I \rightarrow \mathbb{R}^2$ be a smooth closed curve.  Then for any term $P_\nu^\mu\left( k\right)$ with $\nu \geq 2$ that contains derivatives of $k$ of order at most $\ell-1$,
$$\int_I \left| P_\nu^\mu\left( k \right)\right| ds \leq c \, L^{1-\mu-\nu} \left\| k \right\|_2^{\nu-p} \left\| k \right\|_{\ell, 2}^{p}$$
where $p = \frac{1}{\ell} \left( \mu + \frac{1}{2} \nu - 1\right)$ and $c=c\left( \ell, \mu, \nu \right)$.  Moreover, if $\mu+ \frac{1}{2} \nu < 2\ell+1$ then $p<2$ and for any $\varepsilon>0$,
\begin{equation*}
  \int_I  \left| P_\nu^\mu\left( k \right)\right| ds \leq \varepsilon \int_I \left| \partial_s^\ell k \right|^2 ds
   + c\, \varepsilon^{\frac{-p}{2-p}} \left( \int_I \left| k\right|^2 ds\right)^{\frac{\nu-p}{2-p}} + c\left( \int_I \left| k\right|^2 ds \right)^{\mu+\nu-1}\mbox{.}
\end{equation*}
\end{interp}

We conclude this section with the evolution equations for some geometric quantities under the normal curvature flow
\begin{equation} \label{E:generalflow}
  \frac{\partial \gamma}{\partial t} = - F\, \nu \mbox{.}
  \end{equation}
Here $\nu$ is a smooth choice of unit normal vector and the sign is chosen to ensure that \eqref{E:generalflow} is parabolic in the generalised sense.  Throughout this article $L=L\left[ \gamma\right]$ will denote the length of the curve $\gamma$.

The following evolution equations are straightforward to derive using techniques as in \cite{WW}, for example.

\begin{evlneqns} \label{T:evlneqns}
Under the flow \eqref{E:theflow} we have the following evolution equations:
\begin{enumerate}
  \item[\textnormal{(i)}] $\frac{d}{d t} L = -\int_\gamma k \, F\, ds$;\\
  
\noindent For each $\ell = 0, 1, 2, \ldots$, 
  \item[\textnormal{(ii)}] $\frac{\partial}{\partial t} k_{s^\ell} = F_{s^{\ell+2}} + \sum_{j=0}^\ell \partial_{s^j} \left( k\, k_{s^{\ell-j}}F\right)$.
\end{enumerate}  
\end{evlneqns}

\section{The polyharmonic curve diffusion flow} \label{S:PHCD}
In this section we establish our result for the polyharmonic curve diffusion flow \eqref{E:PHCD} for each fixed $m\in \mathbb{N}\cup \left\{ 0 \right\}$.  When $m=0$ we have the classical curve diffusion flow.  The case $m=1$ can be considered the geometric triharmonic heat flow of curves in view of the relationship between curvature and derivatives of $\gamma$.  Since curvature depends upon second spatial derivatives of $\gamma$, the flow \eqref{E:PHCD} has order $2m+4$.

\newtheorem{PHCDL}{Lemma}[section]
\newtheorem{PHCDexp}[PHCDL]{Proposition}
\newtheorem{PHCDL2}[PHCDL]{Corollary}
\newtheorem{PHCDL3}[PHCDL]{Corollary}

\begin{PHCDL} \label{T:PHCDL}
While a solution to the flow \eqref{E:PHCD} with generalised Neumann boundary conditions exists, we have
$$\frac{d}{dt} L\left( t\right) = - \int_\gamma k_{s^{m+1}}^2 ds$$
\end{PHCDL}
\noindent \textbf{Proof:} The result follows directly using Lemma \ref{T:evlneqns} (i) and $m+1$ integrations by parts, noting in each case the boundary term will contain an odd derivative of $k$ that is equal to zero in view of the boundary conditions.\hspace*{\fill}$\Box$\\[8pt]

\noindent \textbf{Remark:} In fact, under the flow \eqref{E:PHCD}, the length $L\left( t\right)$ is strictly decreasing unless $\gamma$ is a straight line segment, because the only smooth solutions that satisfy $k_{s^{m+1}} \equiv 0$ and the boundary conditions are horizontal line segments.\\[8pt]

In view of Lemma \ref{T:PHCDL} and the separation $d_0$ of the supporting parallel lines $\eta_{\pm}$, the length $L\left( t\right)$ of the evolving curve $\gamma\left( \cdot, t\right)$ remains bounded above and below under the flow \eqref{E:PHCD}.\\  

Next we show directly that, provided initially small, $\int k^2 ds$ decays exponentially under the flow.  As in the statement of the main theorem, $\kappa$ denotes the curvature of the initial curve $\gamma_0$.

\begin{PHCDexp} \label{T:PHCDexp}
There exists a constant $\varepsilon>0$, depending only on $m$, such that, if $\gamma_0$ satisfies
\begin{equation} \label{E:PHCDsmallness}
  \int \kappa^2 ds \leq \varepsilon \mbox{,}
\end{equation}
then, while a solution to \eqref{E:PHCD} exists, 
$$ \int k^2 ds  \leq \int \kappa^2 ds \cdot \exp\left( -\delta t \right) \mbox{.}$$
Here $\delta>0$ depends on $\varepsilon$ and $L_0$, the length of $\gamma_0$.
\end{PHCDexp}

\noindent \textbf{Proof:} Under the flow \eqref{E:PHCD}, a straightforward computation using Lemma \ref{T:evlneqns}, integration by parts and the boundary conditions shows that
\begin{equation} \label{E:k2}
  \frac{d}{dt} \int k^2 ds =-2 \int k_{s^{m+2}}^2 ds + \int \left( k^3 \right)_{s^{m+1}} k_{s^{m+1}} ds
   = -2 \int k_{s^{m+2}}^2 ds + \int P_4^{2m  + 2}\left( k\right) ds \mbox{.}
\end{equation}
Since the highest order derivative in $P_4^{2m+2}\left( k\right)$ is $k_{s^{m+1}}$, we have using Proposition \ref{T:int}
$$\int P_4^{2m+2}\left( k \right) ds \leq c\, L^{-\left( 2m+5\right)} \left\| k \right\|_2^{\frac{2m+5}{m+2}} \left\| k \right\|_{m+2, 2}^{\frac{2m+3}{m+2}} \mbox{.}$$
In view of Lemma \ref{T:PSW},
$$\left\| k \right\|_{m+2, 2} \leq c\left( m\right) L^{m + \frac{5}{2}} \left( \int k_{s^{m+2}}^2 ds \right)^{\frac{1}{2}} \mbox{.}$$
We now estimate, again using Lemma \ref{T:PSW},
$$\left\| k \right\|_2^{\frac{2m+5}{m+2}} = L^{\frac{2m+5}{2m+4}} \int k^2 ds \left( \int k^2 ds \right)^{\frac{1}{2m+4}} \leq \frac{1}{\pi} L^{2+ \frac{1}{2m+4}} \int k^2 ds \left( \int k_{s^{m+2}}^2 ds \right)^{\frac{1}{2m+4}} \mbox{.}$$
Combining these, we have 
\begin{equation} \label{E:P4}
  \int P_4^{2m+2}\left( k \right) ds \leq c\, L \int k^2 ds \int k_{s^{m+2}}^2 ds
  \end{equation}
and so from \eqref{E:k2},
$$\frac{d}{dt} \int k^2 ds \leq \left( -2 + c\, L \int k^2 ds \right) \int k_{s^{m+2}}^2 ds \mbox{.}$$
Suppose initially $c\, L \int k^2 ds \leq 2- 2 \tilde \delta$, for some $\tilde \delta >0$.  Then, at least for a short time, $c\, L \int k^2 ds \leq 2- \tilde \delta$.  While this is the case,
$$\frac{d}{dt} \int k^2 ds \leq - \tilde \delta \int k_{s^{m+2}}^2 ds \leq - \tilde \delta \left( \frac{\pi^2}{L^2} \right)^{m+2} \int k^2 ds \leq - \tilde \delta \left( \frac{\pi^2}{L_0^2} \right)^{m+2} \int k^2 ds$$
where we have used again Lemma \ref{T:PSW} and also Lemma \ref{T:PHCDL}, with $L_0$ denoting the length of $\gamma_0$.  The result follows.\hspace*{\fill}$\Box$\\[8pt]

\noindent \textbf{Remarks:} 
\begin{enumerate}
  \item Without the smallness requirement, we can show similarly as in \cite[Theorem 3.1]{DKS02} that if the maximal existence time $T$ of a solution to \eqref{E:PHCD} is finite, then the curvature must blow up in $L^2$.  Specifically, using the second statement of Proposition \ref{T:int} we have from \eqref{E:k2}
  $$ \frac{d}{dt} \int k^2 ds \leq c\, \left( \int k^2 ds\right)^{2m+5}$$
  from which it follows that
  $$\int k^2 ds \geq c\left( T-t\right)^{-\frac{1}{2m+4}} \mbox{.}$$  
  
\item The smallness requirement here  may be compared with the requirement for exponential convergence in the case of closed curves evolving by the polyharmonic curvature flow \cite{PW16}.  Denoting by $\overline{k}$ the average of the curvature over a closed curve, the scale-invariant $\left\| k \right\|_2$ is replaced by
$$K_{osc} = L \int \left( k - \overline{k} \right)^2 ds$$
and a smallness condition on this quantity (together with a condition on the isoperimetric ratio) facilitates exponential convergence of $K_{osc}$.  The two quantities again appear in parallel in the Poincar\'{e}-Sobolev-Wirtinger inequalities for curves with boundary and closed curves.\\[8pt]
\end{enumerate}

The next step is to show that under the flow, $L^2$ norms of all curvature derivatives remain bounded.  This proof here is considerably more direct than in the subsequent section for the flow \eqref{E:theflow}.

\begin{PHCDL2} \label{T:PHCDL2}
Suppose $\gamma_0$ satisfies the conditions of Theorem \ref{T:main1} including the smallness condition \eqref{E:PHCDsmallness}.  Then, while a solution to the flow \eqref{E:PHCD} with generalised Neumann boundary conditions exists, we have for all $\ell \in \mathbb{N}\cup \left\{ 0\right\}$,
$$\int k_{s^\ell}^2 ds \leq C_{m, \ell} \mbox{,}$$
for constants $C_{m, \ell}$.
\end{PHCDL2}

\noindent \textbf{Proof:} Under the flow \eqref{E:PHCD}, a straightforward computation using Lemma \ref{T:evlneqns} and integration by parts with repeated application of the boundary conditions and the consequence \eqref{E:bdy}, gives that for each $\ell$,
\begin{equation} \label{E:evlnL2kl}
  \frac{d}{dt} \int k_{s^\ell}^2 ds = -2 \int k_{s^{m+\ell+2}}^2 ds + \int P_4^{2m + 2\ell + 2}\left( k\right) ds \mbox{,}
\end{equation}
where the highest order of derivatives of $k$ in the second above term is $m+\ell+1$.  Using Proposition \ref{T:int} we have for any $\varepsilon >0$,
\begin{equation} \label{E:P4}
  \int P_4^{2m + 2\ell + 2}\left( k\right) ds \leq \varepsilon \int k_{s^{m+\ell+2}}^2 ds + c\left( m, \ell, \varepsilon \right) \left( \int k^2 ds\right)^{2m + 2\ell +5} \mbox{,}
\end{equation}
so from \eqref{E:evlnL2kl} we obtain by taking $\varepsilon=1$
$$ \frac{d}{dt} \int k_{s^\ell}^2 ds \leq - \int k_{s^{m+\ell+2}}^2 ds + c \left( \int k^2 ds\right)^{2m + 2\ell +5} \mbox{.}$$
Using now Lemma \ref{T:PSW} and Lemma \ref{T:PHCDL} we obtain
$$ \frac{d}{dt} \int k_{s^\ell}^2 ds \leq - \left( \frac{\pi^2}{L_0^2}\right)^{m+2} \int k_{s^\ell}^2 ds + c \left( \int k^2 ds\right)^{2m + 2\ell +5}$$
from which the result follows since $\int k^2 ds$ is bounded in view of Proposition \ref{T:PHCDexp}.\hspace*{\fill}$\Box$\\[8pt]

\noindent \textbf{Remark:} In view of Corollary \ref{T:PHCDL2}, a standard contradiction argument using short-time existence implies that in fact the solution to \eqref{E:PHCD} exists for all time, that is, $T=\infty$.\\[8pt]

Proposition \ref{T:PHCDexp} and Corollary \ref{T:PHCDL2} imply via interpolation that all curvature derivatives decay exponentially in $L^2$ and, via Lemma \ref{T:PSW}, in $L^\infty$.

\begin{PHCDL3} \label{T:PHCDL3}
Suppose $\gamma_0$ satisfies the conditions of Theorem \ref{T:main1} including the smallness condition \eqref{E:PHCDsmallness}.  Under the flow \eqref{E:PHCD}, there exist $\delta_{\ell, m}>0$, depending only on $\varepsilon$ and $L_0$ such that, for all $\ell \in \mathbb{N} \cup \left\{ 0 \right\}$,
$$\int k_{s^\ell}^2 ds \leq \tilde C_{m, \ell} \exp \left( - \delta_{\ell, m} t\right) \mbox{.}$$
The quantities $\left\| k_{s^\ell} \right\|_\infty$ also decay exponentially for all $\ell$.
\end{PHCDL3}

\noindent \textbf{Proof:} The curvature derivative decay in $L^2$ follows by standard integration by parts; we give the first two calculations:
$$\int k_s^2 ds = \left[ k_s k \right]_{\partial \gamma} - \int k_{ss} k \, ds \leq \left( \int k_{ss}^2 ds\right)^{\frac{1}{2}} \left( \int k^2 ds\right)^{\frac{1}{2}} \mbox{.}$$
Here we have used that the boundary term contains an odd derivative of $k$ so is equal to zero.  The $\int k_{ss}^2 ds$ factor is bounded by Corollary \ref{T:PHCDL2} and Lemma \ref{T:PHCDL} so then the exponential convergence follows from Proposition \ref{T:PHCDexp}.

We next compute
$$\int k_{ss}^2 ds = \left[ k_{ss} k_s \right]_{\partial \gamma} - \int k_{sss} k_s \, ds \leq \left( \int k_{sss}^2 ds\right)^{\frac{1}{2}} \left( \int k_s^2 ds\right)^{\frac{1}{2}} \mbox{.}$$
Again we have used that the boundary term contains an odd derivative of $k$.  The $\int k_{sss}^2 ds $ factor is bounded by Corollary \ref{T:PHCDL2} and Lemma \ref{T:PHCDL} so then the exponential convergence follows from the previous step.

We can continue this way to obtain all curvature derivatives in $L^2$ decay exponentially.  The $L^\infty$ decay then follows from Lemma \ref{T:PSW} and Lemma \ref{T:PHCDL}.\hspace*{\fill}$\Box$
\mbox{}\\[8pt]

Using now the evolution equation \eqref{E:PHCD} we obtain uniform bounds on all derivatives of the immersion $\gamma: \left[ -1, 1\right] \times \left[ 0, \infty\right) \rightarrow \mathbb{R}^2$.  This implies there exists an immersion $\gamma_\infty : \left[ -1, 1\right] \rightarrow \mathbb{R}^2$ satisfying the boundary conditions and a subsequence $t_j \rightarrow \infty$ such that $\gamma\left( \cdot, t_j \right) \rightarrow \gamma_\infty$ in $C^\infty\left( \left[ -1, 1\right], \mathbb{R}^2 \right)$.  Since $\left\| k \right\|_\infty \rightarrow 0$, the curve $\gamma_\infty$ is a straight line segment which is horizontal in view of the Neumann boundary condition.  Exponential convergence in $C^\infty$ of $\gamma$ to $\gamma_\infty$ now follows by the same argument as in \cite{MWW} using exponential convergence of the curvature and its derivatives.  This completes the proof of Theorem \ref{T:main1}. \hspace*{\fill}$\Box$\\[8pt]
 
 \noindent \textbf{Remark:} While we don't know the precise height of the limiting straight horizontal line segment, that $\left\| k_{s^{2m+2}} \right\|_\infty$ decays exponentially shows that the solution curve remains within a bounded distance of the initial curve: for any $x$,
$$\left| \gamma\left( x, \tilde t\right) - \gamma\left( x, 0 \right) \right| \leq \int_0^{\tilde t} \left| \frac{\partial \gamma}{\partial t} \left( x, t\right) \right| dt \leq c \int_0^{\tilde t} e^{-\delta \, t} dt = \frac{c}{\delta} \left( 1- e^{-\delta \tilde t} \right) \mbox{.}$$
 
\section{The gradient flow for $\int k_{s^m}^2 ds$} \label{S:setup}
\newtheorem{BCs}{Lemma}[section]
\newtheorem{BD2}[BCs]{Lemma}
\newtheorem{Length}[BCs]{Lemma}
\newtheorem{L2bounds}[BCs]{Proposition}
\newtheorem{kexp}[BCs]{Proposition}
\newtheorem{dexp}[BCs]{Corollary}
For a suitably smooth curve $\gamma$ as described in Section \ref{S:intro}, we are interested in the associated curvature-dependent energies
$$E\left[ \gamma\right] = \frac{1}{2} \int_\gamma k_{s^m}^2\, ds$$
and the corresponding $L^2$-gradient flows with suitable associated generalised Neumann boundary conditions.  As the energy involves $m$ derivative of curvature, so $m+2$ derivatives of $\gamma$, the gradient flow will be of order $2m+4$.

Under a normal variation $\tilde \gamma = \gamma + \varepsilon F \nu$
straightforward calculations yield
\begin{multline} \label{E:1}
  \left. \frac{d}{d\varepsilon} E\left[ \tilde \gamma \right] \right|_{\varepsilon=0}
  = - 2 \int_\gamma \left[ \left( -1\right)^{m+1} k_{s^{2m+2}} - \sum_{j=1}^m \left( -1\right)^j k\, k_{s^{m+j}} k_{s^{m-j}} - \frac{1}{2} k\, k_{s^m}^2 \right] F \, ds\\
  + 2 \left[ \sum_{j=0}^{m+1} \left( -1\right)^j k_{s^{m+j}} \partial_{s^{m+1-j}} F + \sum_{j=1}^{m} \sum_{\ell=0}^{j-1} \left( -1\right)^\ell k_{s^{m+\ell}} \partial_{s^{j-1-\ell}} \left( k\, k_{s^{m-j}} F \right) \right]_{\partial \gamma} \mbox{.}
  \end{multline}
In particular, the above follows from the variations
$$\left. \frac{\partial}{\partial \varepsilon} k_{s^m} \right|_{\varepsilon=0}= \partial_{s^{m+2}} F + \sum_{j=0}^m \partial_{s^j} \left( k\, k_{s^{m-j}} F \right)$$
and 
$$\left. \frac{\partial}{\partial \varepsilon} ds \right|_{\varepsilon=0} = -k\, F \, ds \mbox{,}$$
where the boundary terms in \eqref{E:1} appear via repeated integration by parts.  Details behind these calculations may be found for example in \cite{WW}.

`Natural boundary conditions' for the corresponding $L^2$-gradient flow would
ensure that the boundary term in \eqref{E:1} is equal to zero.  Assuming for now it is equal to zero we would take the normal flow speed
\begin{equation} \label{E:speed}
  F=   \left( -1\right)^{m+1} k_{s^{2m+2}} - \sum_{j=1}^m \left( -1\right)^j k\, k_{s^{m+j}} k_{s^{m-j}} - \frac{1}{2} k\, k_{s^m}^2
  \end{equation}
and the corresponding $L^2$-gradient flow is then \eqref{E:theflow}.\\

Let us now establish a mathematically-reasonably choice of boundary conditions.  Beginning with the classical Neumann boundary condition and differentiating in time (see also \cite[Lemma 2.5]{WW} for example) we have 
\begin{equation} \label{E:Fs}
  F_s\left( \pm 1, t\right) = 0 \mbox{.}
\end{equation}

If we assume as in previous work the `no curvature flux condition' at the boundary, 
\begin{equation} \label{E:noflux}
  k_s\left( \pm 1, t\right) =0 \mbox{,}
\end{equation}  
then from the evolution equation for $k_s$,
$$\frac{\partial}{\partial t} k_s = F_{s^3} + \sum_{j=0}^1 \partial_{s^j} \left( k\, k_{s^{1-j}} F \right)$$
we see that, on the boundary, we must also have $F_{s^3}\left( \pm 1, t\right) \equiv 0$.  More generally, by similar arguments in turn, for each odd derivative of $k$ equal to zero on the boundary, we see that the next odd derivative of $F$ is also equal to zero on the boundary.  Assuming then that all odd derivatives of $k$ up to order $2m+1$ are equal to zero on the boundary and taking into account the corresponding behaviour of the odd derivatives of $F$ on the boundary, we see that this choice of boundary conditions does render the boundary term in \eqref{E:1} equal to zero.  Moreover, we then have by an inductive argument similar to that in \cite{WW}:

\begin{BCs} \label{T:BCs}
With classical Neumann conditions and all odd derivatives of $k$ up to order $2m+1$ equal to zero on the boundary, a solution to the flow \eqref{E:theflow} satisfies $k_{s^{2\ell-1}}=0$ and $F_{s^{2\ell -1}} =0$ on the boundary for all $\ell \in \mathbb{N}$.
\end{BCs}

Throughout our work it will be necessary to check that various boundary terms arising by parts are equal to zero.  These boundary terms typically are sums of products.  In many cases, each product has three factors that are each either an \emph{even} or an \emph{odd} number of iterated spatial derivatives of $k$ or of $F$.  For certain products, we need to know that an odd number of derivatives always produces an odd iterated derivative factor of $k$ or of $F$, while for others we need to know that an even number of derivatives always produces an odd derivative factor.  The presence of such a factor in each term then ensures by Lemma \ref{T:BCs} that the boundary term is equal to zero.

We introduce the notation $e$ to denote an even (or zeroth order) derivative factor of $k$ or of $F$ and $o$ to denote an odd derivative factor of $k$ or of $F$.  We allow sums of such factors in the notation $e\, o$ etc.  We also use subscripts $e$ and $o$ to denote respectively an even (or zero) number or an odd number of spatial derivatives.  Some of the results that we will need later can now be stated as follows:

\begin{BD2} \label{T:BD2}
Terms of the form $\left( e\, e\, e\right)_o$, $\left( e\, o\, o\right)_o$ and $\left( e\, e\, o \right)_e$ always contain an $o$ factor.
\end{BD2}

\noindent \textbf{Proof:} Using the product rule we begin with
$$\left( e\, e\, e\right)_s = e\, e\, o \mbox{,}$$
which has the required form.  Differentiating a second time,
$$\left( e\, e\, e\right)_{ss} = \left( e\, e\, o \right)_s = e\, o\, o  +   e\, e\, e \mbox{,}$$
and a third time
$$\left( e\, e\, e\right)_{s^3}=  o\, o\, o  +  e\, e\, o  \mbox{.}$$
We see that the third derivative also consists of terms with an $o$ factor.  Continuing
$$\left( e\, e\, e\right)_{s^4}=  o\, o\, e  +  e\, o\, o  +  e\, e\, e $$
and
$$\left( e\, e\, e\right)_{s^5}=  o\, o\, o  +  o\, e\, e  \mbox{,}$$
which is the same form as the third derivative.  We conclude the result for all odd derivatives $\left( e\, e\, e\right)_o$.\\

For $\left( e\, o\, o\right)_o$, we begin with
$$\left( e\, o\, o\right)_s = e\, e\, o + o\, o\, o\mbox{,}$$
which has the required form.  Differentiating a second time,
$$\left( e\, o\, o\right)_{ss} = e\, o\, o + e\, e\, e$$
and a third time
$$\left( e\, o\, o\right)_{s^3} = e\, e\, o + o\, o\, o \mbox{,}$$
which is the same form as the first derivative, so we conclude the result for all odd derivatives $\left( e\, o\, o\right)_o$.\\

For $\left( e\, e\, o\right)$, the zeroth order derivative has an $o$ factor as required.  We compute
$$\left( e\, e\, o\right)_s =  e\, o\, o  + e\, e\, e $$
and so
$$\left( e\, e\, o\right)_{ss} = e\, e\, o  + o\, o\, o \mbox{;}$$
so each term has an $o$ factor as required.  Continuing
$$\left( e\, e\, o\right)_{s^3} =  e\, o\, o + e\, e\, e \mbox{;}$$
this is the same form as $\left( e\, e\, o\right)_s$, so therefore even derivatives, like $\left( e\, e\, o\right)_{ss}$, will always consist of sums of terms containing $o$ factors, as required.\hspace*{\fill}$\Box$ 
\mbox{}\\[8pt]

\noindent \textbf{Remark:} We will occasionally need results related to the above, in particular when `square' factors appear in the products to be considered.  We will develop those results directly where they are needed to follow.\\

Notwithstanding our earlier comments on short-time existence, our first result for the flow \eqref{E:theflow} shows that if the initial energy is small, then the length of the evolving curve does not increase.

\begin{Length} \label{T:Length}
If the initial curve $\gamma_0$ has sufficiently small energy \eqref{E:E} depending only on $m$, then, under the flow \eqref{E:theflow} with normal speed \eqref{E:speed}, the length of $\gamma$ does not increase.
\end{Length}

\noindent \textbf{Proof:} We have using Lemma \ref{T:evlneqns}, (i),
\begin{equation} \label{E:Lksm}
  \frac{d}{dt} L = - \int k \left[ \left( -1\right)^{m+1} k_{s^{2m+2}} + \sum_{j=1}^m \left( -1\right)^{j+1} k\, k_{s^{m+j}} k_{s^{m-j}} - \frac{1}{2} k\, k_{s^m}^2 \right] ds \mbox{.}
\end{equation}
By integrating by parts the first term $m+1$ times and using the boundary conditions, we have 
$$-\left( -1\right)^{m+1}\int k  \, k_{s^{2m+2}} ds = - \int k_{s^{m+1}}^2 ds \mbox{,}$$
while integrating by parts the $j$th term in the sum $j$ times we can see that the rest of the terms have the form $\int P_4^{2m}\left( k\right) ds$ with the highest order derivative of $k$ being $k_{s^m}$.  Thus we estimate using Proposition \ref{T:int}
$$\int P_4^{2m}\left( k\right) ds \leq c\left( m\right) L^{-\left( 2m+3\right)} \left\| k \right\|_2^{\frac{2m+3}{m+1}} \left\| k\right\|_{m+1, 2}^{\frac{2m+1}{m+1}} \mbox{.}$$
Using now Lemma \ref{T:PSW} we have
$$ \left\| k\right\|_{m+1, 2} \leq c\left( m\right) L^{m+ \frac{3}{2}} \left( \int k_{s^{m+1}}^2 ds \right)^{\frac{1}{2}}$$
and
$$\int k^2 ds \leq \left( \frac{L^2}{\pi^2}\right)^{m+1}  \left( \int k_{s^{m+1}}^2 ds \right) \mbox{,}$$
so
$$\int P_4^{2m}\left( k\right) ds \leq c\left( m\right) \left\| k\right\|_2^2 \int k_{s^{m+1}}^2 ds$$
and from \eqref{E:Lksm} we obtain
$$ \frac{d}{dt} L \leq \left( -1 + c\left( m\right)  L^{2m+1} \int k_{s^{m}}^2 ds \right) \int k_{s^{m+1}}^2 ds \mbox{.}$$
Since $\int k_{s^{m}}^2 ds$ is nonincreasing under the flow by construction, it follows that if $\gamma_0$ has\\ $L^{2m+3} \int k_{s^{m}}^2 ds$ sufficiently small, depending only on $m$, then $L$ does not increase under the flow \eqref{E:theflow}. \hspace*{\fill}$\Box$

\begin{L2bounds} \label{T:L2bounds}
If the initial curve $\gamma_0$ has sufficiently small energy \eqref{E:E} depending only on $m$, then, under the flow \eqref{E:theflow} with normal speed \eqref{E:speed}, there are constants $C_{m, \ell}$ depending only on $L_0$ and the initial energy such that
$$\int k_{s^\ell}^2 ds \leq C_{m, \ell}$$
for all $\ell \in \mathbb{N}\cup \left\{ 0\right\}$.
\end{L2bounds}

\noindent \textbf{Proof:} For $\ell \leq m$, the result is immediate via Lemma \ref{T:PSW} and Lemma \ref{T:Length}.  For any $\ell$, we have via Lemma \ref{T:evlneqns}
\begin{equation} \label{E:ksl2}
  \frac{d}{dt} \int k_{s^\ell}^2 ds = 2 \int k_{s^\ell} \left[ F_{s^{\ell+2}} + \sum_{j=0}^\ell \partial_{s^j} \left( k\, k_{s^{\ell-j}} F\right) \right] ds -  \int k_{s^\ell}^2 k\, F \, ds \mbox{.}
  \end{equation}
For each $\ell >m$ we examine each of the terms on the right hand side of \eqref{E:ksl2} in turn.  Since the leading term of $F_{s^{\ell+2}}$ is $k_{s^{2m+\ell+4}}$, we will integrate by parts the first term in \eqref{E:ksl2} $\left( m+2\right)$ times:
$$\int k_{s^\ell} F_{s^{\ell+2}} ds = \left[ k_{s^\ell} F_{s^{\ell+1}} \right]_{\partial \gamma} - \int k_{s^{\ell+1}} F_{s^{\ell+1}} ds \mbox{.}$$
Regardless of $\ell$, the boundary term above will have an odd derivative, so is equal to zero by Lemma \ref{T:BCs}.  Integrating by parts again,
$$\int k_{s^\ell} F_{s^{\ell+2}} ds = \left[ k_{s^{\ell+1}} F_{s^{\ell}} \right]_{\partial \gamma} - \int k_{s^{\ell+2}} F_{s^{\ell}} ds \mbox{.}$$
Again, the boundary term will have an odd derivative so is equal to zero.  With a further $m$ integrations by parts, observing that the boundary terms are always equal to zero, we obtain
\begin{align} \label{E:ksl2aux} \nonumber
 & \int k_{s^\ell} F_{s^{\ell+2}} ds\\ \nonumber
 &= \left( -1\right)^m \int k_{s^{\ell+m+2}} F_{s^{\ell-m}} ds\\ \nonumber
 &= \left( -1\right)^m \int k_{s^{\ell+m+2}} \left[ \left( -1\right)^{m+1} k_{s^{2m+2}} - \sum_{j=1}^m k\, k_{s^{m+j}} k_{s^{m-j}} - \frac{1}{2} k\, k_{s^m}^2 \right]_{s^{\ell - m}}ds\\ \nonumber
  &= -\int k_{s^{\ell+m+2}}^2 ds - \left( -1\right)^m \sum_{j=1}^m \int k_{s^{\ell+m+2}}\left( k\, k_{s^{m+j}} k_{s^{m-j}} \right)_{s^{\ell-m}} ds\\ 
  & \quad- \frac{ \left( -1\right)^m }{2}\int k_{s^{\ell+m+2}} \left( k\, k_{s^m}^2 \right)_{s^{\ell - m}} ds \mbox{.}
  \end{align}
We want to confirm that the summation and last terms here have the form $\int P_4^{2m+2\ell  +2}\left( k\right) ds$ with highest order derivative $k_{s^{\ell+m+1}}$.  Integrating by parts the last term,
\begin{equation} \label{E:extra1}
  \int k_{s^{\ell + m + 2}} \left( k\, k_{s^m}^2 \right)_{s^{\ell-m}} ds = \left[ \left( k\, k_{s^m}^2\right)_{s^{\ell-m}} k_{s^{\ell+m+1}}\right]_{\partial \gamma} - \int k_{s^{\ell+m+1}} \left( k\, k_{s^m}^2 \right)_{s^{\ell -m+1}} ds \mbox{.}
  \end{equation}
For the above boundary term, if $\ell+m+1$ is odd then we have an odd derivative of $k$ factor which is zero by Lemma \ref{T:BCs}.  If, on the other hand, $\ell+m+1$ is even, then $\ell+m$ is odd and so is $\ell-m$.  Using the notation of Lemma \ref{T:BD2}, we have multiplying $k_{s^{\ell-m}}$ a term either of the form $\left( e\, e\, e\right)_o$ or $\left( e\, o\, o\right)_o$.  Lemma \ref{T:BD2} gives that either of these consist of terms all with odd derivatives of $k$, thus the boundary term is equal to zero in this case also.  The remaining integral term in \eqref{E:extra1} has the form $\int P_4^{2m+2\ell+2}\left( k\right) ds$ with highest order derivative $k_{s^{m+\ell+1}}$. 

For the summation term in \eqref{E:ksl2aux} we again need one integration by parts: for each $j$,
\begin{multline*}
  \int k_{s^{\ell+m+2}}\left( k\, k_{s^{m+j}} k_{s^{m-j}} \right)_{s^{\ell-m}} ds\\
   = \left[ k_{s^{\ell+m+1}}\left( k\, k_{s^{m+j}} k_{s^{m-j}} \right)_{s^{\ell-m}} \right]_{\partial \gamma} - \int k_{s^{\ell+m+1}}\left( k\, k_{s^{m+j}} k_{s^{m-j}} \right)_{s^{\ell-m+1}} ds \mbox{.}
 \end{multline*}
If $\ell+m+1$ is odd then again clearly the boundary term is equal to zero.  Otherwise, $\ell - m$ is odd and the other factor in the boundary term has the form $\left( e\, e\, e\right)_o$ or $\left( e\, o\, o\right)_o$.  As before, the boundary term is again equal to zero for each $i$ and the remaining integral terms have the correct form.\\

Returning now to \eqref{E:ksl2}, using \eqref{E:speed} we have
\begin{equation} \label{E:ksl2aux2}
  \int k_{s^\ell}^2 k\, F\, ds = \int k\, k_{s^\ell}^2 \left[ \left( -1\right)^{m+1} k_{s^{2m+2}} - \sum_{j=1}^m k\, k_{s^{m+j}} k_{s^{m-j}} - \frac{1}{2} k\, k_{s^m}^2 \right] ds \mbox{.}
  \end{equation}
  We want to show that all these terms have either the form $\int P_{4}^{2\ell + 2m +2}\left( k \right) ds$ with highest derivative $k_{s^{\ell+m+1}}$ or $\int P_6^{2m+2\ell}\left( k\right) ds$ with no higher derivative than $k_{s^{\ell+m}}$.  The last term above already fits this latter form.  On the first term we will need to integrate by parts $m+1$ times: first
  $$\int k\, k_{s^\ell}^2 k_{s^{2m+2}} ds = \left[ k\, k_{s^\ell}^2 k_{s^{2m+1}} \right]_{\partial \gamma} - \int \left( k\, k_{s^\ell}^2 \right)_s k_{s^{2m+1}} ds \mbox{.}$$
  The boundary term here has an odd derivative of $k$ so is equal to zero.  Integrating by parts a second time,
   $$\int k\, k_{s^\ell}^2 k_{s^{2m+2}} ds = -\left[ \left( k\, k_{s^\ell}^2 \right)_s k_{s^{2m}} \right]_{\partial \gamma} + \int \left( k\, k_{s^\ell}^2 \right)_{ss} k_{s^{2m}} ds \mbox{.}$$
   This time $k_{s^{2m}}$ in the boundary term is an even derivative so we look at the other factor.  Depending on $\ell$, the other factor has the form $\left( e\, e\, e\right)_o$ or $\left( e\, o\, o\right)_o$; so in both cases these will be equal to zero by Lemma \ref{T:BD2} and Lemma \ref{T:BCs}.  With each of the further integrations by parts, the boundary terms will be one of the two types above, so we are left with
    $$\int k\, k_{s^\ell}^2 k_{s^{2m+2}} ds = \left( -1\right)^{m+1} \int \left( k\, k_{s^\ell}^2 \right)_{s^{m+1}} k_{s^{m+1}} ds  = \int P_4^{2m+2\ell+2}\left( k \right) ds\mbox{.}$$
   with no higher derivative of $k$ than $k_{s^{m+\ell+1}}$ appearing.
   
   The terms in the summation in \eqref{E:ksl2aux2} require some integration by parts if $j$ is large relative to $\ell$.  Specifically, for the $j$th term we should do $j$ integrations by parts to ensure no derivative of order higher than $m+\ell$ appears.  We have
   $$\int k^2 k_{s^\ell}^2 k_{s^{m-j}} k_{s^{m+j}} ds = \left[ k^2 k_{s^\ell}^2 k_{s^{m-j}} k_{s^{m+j-1}} \right]_{\partial \gamma} -\int \left( k^2 k_{s^\ell}^2 k_{s^{m-j}} \right)_s k_{s^{m+j-1}} ds  \mbox{.}$$
   The boundary term above will always have an odd derivative of $k$ so is equal to zero.  Continuing in the case $j>1$ we have
   $$\int k^2 k_{s^\ell}^2 k_{s^{m-j}} k_{s^{m+j}} ds =- \left[ \left( k^2 k_{s^\ell}^2 k_{s^{m-j}}\right)_s k_{s^{m+j-2}} \right]_{\partial \gamma} +\int \left( k^2 k_{s^\ell}^2 k_{s^{m-j}} \right)_{ss} k_{s^{m+j-2}} ds  \mbox{.}$$
   If $m+ j$ is odd, then the boundary term is equal to zero.  In the case $m+j$ is even, then $m-j$ is also even and applying the product rule to expand the first derivative, we see that every term will have an odd derivative so is equal to zero.  Clearly, similar terms will arise with further integrations by parts.  We need to see that an odd derivative of a product of the form $\mbox{square } \times \mbox{square } \times e$ is always equal to zero.  Expanding out such a derivative using the binomial theorem, the terms with an odd derivative of $e$ are already equal to zero via Lemma \ref{T:BCs} so we only need to check those terms with an odd derivative of a square.  We have
   $$\left( e\, e\right)_s = e\, o$$
   $$\left( e\, e\right)_{ss} = e\, e + o\, o$$
   $$\left( e\, e\right)_{s^3} = e\, o$$
   and so generally $\left( e\, e\right)_o = e\, o$.  Similarly
   $$\left( o\, o\right)_s = e\, o$$
   so the same pattern will occur and we see that all odd derivatives of squares contain an $o$ factor, thus the boundary terms generated through repeated integration by parts will always be equal to zero.  Therefore we have
   $$\int k^2 k_{s^\ell}^2 k_{s^{m-j}} k_{s^{m+j}} ds =\left( -1\right)^j \int \left( k^2 k_{s^\ell}^2 k_{s^{m-j}} \right)_{s^j} k_{s^{m}} ds = \int P_6^{2m+2\ell}\left( k\right) ds\mbox{,}$$
   where no higher derivative than $k_{s^{m+\ell}}$ appears on the right hand side.\\
   
   It remains now to consider the summation term in \eqref{E:ksl2}.  For $j=0$ this term has the same form as that dealt with above, so assume $1\leq j \leq \ell$.  To handle the boundary terms arising from integration by parts, it is going to be easiest to do $j$ integrations by parts first without substituting in the form of $F$.  The first integration by parts gives
   $$\int k_{s^\ell} \partial_{s^j} \left( k\, k_{s^{\ell-j}} F\right) ds = \left[ k_{s^\ell} \partial_{s^{j-1}} \left( k\, k_{s^{\ell-j}} F\right) \right]_{\partial \gamma} - \int k_{s^{\ell+1}} \partial_{s^{j-1}} \left( k\, k_{s^{\ell-j}} F\right) ds \mbox{.}$$
If $\ell$ is odd, then the above boundary term is clearly equal to zero.  If $\ell$ is even then we consider two cases.  If $j$ is also even then the other factor in the boundary term has the form $\left( e\, e\, o\right)_o$, that by Lemma \ref{T:BD2} always contains an odd derivative of $k$ or of $F$.  If $k$ is odd then the same factor instead has the form $\left( e\, e\, e\right)_o$ that by Lemma \ref{T:BD2} also always contains an odd derivative of $k$ or of $F$.  We conclude by Lemma \ref{T:BCs} that the boundary term in all cases is equal to zero.  

Integrating by parts a second time (for $j>1$) we have
$$\int k_{s^\ell} \partial_{s^j} \left( k\, k_{s^{\ell-j}} F\right) ds = - \left[ k_{s^\ell+1} \partial_{s^{j-2}} \left( k\, k_{s^{\ell-j}} F\right) \right]_{\partial \gamma} - \int k_{s^{\ell+2}} \partial_{s^{j-2}} \left( k\, k_{s^{\ell-j}} F\right) ds \mbox{.}$$
Here the boundary term is clearly equal to zero if $\ell$ is even, while if $\ell$, two cases depending on $j$ are handled as above.

Continuing like this we have for each $j=1, \ldots, \ell$,
\begin{multline} \label{E:intFin}
  \int k_{s^\ell} \partial_{s^j} \left( k\, k_{s^{\ell-j}} F\right) ds = \left( -1\right)^j \int k_{s^{\ell+j}} k\, k_{s^{\ell-j}} F ds\\
  =  \left( -1\right)^j \int k_{s^{\ell+j}} k\, k_{s^{\ell-j}} \left[ \left( -1\right)^m k_{s^{2m+2}} - \sum_{i=1}^m k\, k_{s^{m+i}} k_{s^{m-i}} - \frac{1}{2} k\, k_{s^m}^2 \right] ds \mbox{.}
  \end{multline}
 Let us now deal with each of these terms separately, remembering that we want factors of the form $ \int P_4^{2m+2\ell+2}\left( k \right) ds$ with no higher derivative than $k_{s^{m+ \ell + 1}}$ and $\int P_6^{2m+2\ell} \left( k \right) ds$ with no higher derivative than $k_{s^{m+\ell}}$.  The last term has the correct form for $j\leq m$ but if $j>m$ needs $m-j$ integrations by parts.  In this case we have
 $$\int k^2 k_{s^m}^2 k_{s^{\ell-j}} k_{s^{\ell+j}} ds = \left[ k^2 k_{s^m}^2 k_{s^{\ell-j}} k_{s^{\ell+j-1}} \right]_{\partial \gamma} - \int k_{s^{\ell + j-1}} \left( k^2 k_{s^m}^2 k_{s^{\ell-j}} \right)_s ds \mbox{;}$$
 the boundary term clearly contains an odd derivative of $k$ so is equal to zero by Lemma \ref{T:BCs}.  Integrating by parts again (if $j>m+1$),
 $$\int k^2 k_{s^m}^2 k_{s^{\ell-j}} k_{s^{\ell+j}} ds = -\left[ \left( k^2 k_{s^m}^2 k_{s^{\ell-j}} \right)_s k_{s^{\ell+j-2}} \right]_{\partial \gamma} - \int k_{s^{\ell + j-2}} \left( k^2 k_{s^m}^2 k_{s^{\ell-j}} \right)_{ss} ds \mbox{.}$$
 Now, if $\ell+j$ is odd then the boundary term is equal to zero via Lemma \ref{T:BCs}.  If $\ell+j$ is even then so is $\ell-j$ and the other factor in the boundary term has the form $\left( \mbox{square } \times \mbox{ square }\times e \right)_o$.  As we saw before, such terms always contain an odd derivative factor, so are equal to zero by Lemma \ref{T:BCs}.  Integrating by parts once again (if $j>m+2$),
  $$\int k^2 k_{s^m}^2 k_{s^{\ell-j}} k_{s^{\ell+j}} ds = \left[ \left( k^2 k_{s^m}^2 k_{s^{\ell-j}} \right)_{ss} k_{s^{\ell+j-3}} \right]_{\partial \gamma} - \int k_{s^{\ell + j-3}} \left( k^2 k_{s^m}^2 k_{s^{\ell-j}} \right)_{s^3} ds \mbox{.}$$
Here, if $\ell+j$ is even then the boundary term is again equal to zero.  If $\ell+j$ is odd then the other factor in the boundary term has the form $\left( \mbox{square } \times \mbox{ square }\times o \right)_e$.  Expanding out such a term using the binomial theorem, all terms with $e$ derivatives of $o$ are fine so we need only consider the terms with $o$ derivatives of $o$.  For such terms one square factor has an $o$ derivative, so, as before, this will generate an $o$ and again the boundary term is equal to zero by Lemma \ref{T:BCs}.  Continuing in this way, we can write for each integer $j\in \left( m, \ell\, \right]$,
$$\int k^2 k_{s^m}^2 k_{s^{\ell-j}} k_{s^{\ell+j}} ds = \left( -1\right)^{j-m} \int \left( k^2 k_{s^m}^2 k_{s^{\ell-j}} \right)_{s^{j-m}} k_{s^{\ell-\left( j-m\right)}} ds = \int P_6^{2m+2\ell}\left( k\right) ds \mbox{,}$$
with no higher derivative appearing than $k_{s^{\ell+m}}$ appearing on the right hand side.

Turning now to the first term of \eqref{E:intFin}, if $j\leq m+1$ then this term already has the form $\int P_4^{2m + 2\ell+2}\left( k\right) ds$ with no higher derivative than $k_{s^{m+\ell +1}}$ appearing.  On the other hand,
if $j>m+1$ then we will need to perform $j-\left( m+1\right)$ integrations by parts.  We have
$$\int k_{s^{\ell+j}} k\, k_{s^{\ell-j}} k_{s^{2m+2}} ds = \left[  k\, k_{s^{\ell-j}} k_{s^{2m+2}} k_{s^{\ell+j-1}} \right]_{\partial \gamma} - \int \left( k\, k_{s^{\ell-j}} k_{s^{2m+2}}\right)_s k_{s^{\ell+j-1}} ds \mbox{.}$$
The boundary term here clearly has an odd derivative of $k$ so is equal to zero by Lemma \ref{T:BCs}.  Integrating by parts again, if necessary,
$$\int k_{s^{\ell+j}} k\, k_{s^{\ell-j}} k_{s^{2m+2}} ds = - \left[ \left( k\, k_{s^{\ell-j}} k_{s^{2m+2}}\right)_s k_{s^{\ell+j-2}} \right]_{\partial \gamma} + \int \left( k\, k_{s^{\ell-j}} k_{s^{2m+2}}\right)_{ss} k_{s^{\ell+j-2}} ds \mbox{.}$$
Here, if $\ell+j$ is odd then the boundary term is clearly equal to zero.  If $\ell+j$ is even, then so is $\ell-j$ and the other factor in the boundary term has the form $\left( e\, e\, e\right)_o$; Lemma \ref{T:BD2} and Lemma \ref{T:BCs} then give that the boundary term is again equal to zero.  Integrating by parts once more, if necessary,
$$\int k_{s^{\ell+j}} k\, k_{s^{\ell-j}} k_{s^{2m+2}} ds =  \left[ \left( k\, k_{s^{\ell-j}} k_{s^{2m+2}}\right)_{ss} k_{s^{\ell+j-3}} \right]_{\partial \gamma} - \int \left( k\, k_{s^{\ell-j}} k_{s^{2m+2}}\right)_{s^3} k_{s^{\ell+j-3}} ds \mbox{.}$$
In this case the boundary term is clearly equal to zero if $\ell+j$ is even, while if $\ell +j $ is odd then the other factor has the form $\left( e\, o\, e\right)_e$; again Lemma \ref{T:BD2} and \ref{T:BCs} give that the boundary term is equal to zero.  Continuing in this way, we have for each $j>m+1$,
$$\int k_{s^{\ell+j}} k\, k_{s^{\ell-j}} k_{s^{2m+2}} ds = \left( -1\right)^{j-\left( m+1\right)} \int \left( k\, k_{s^{\ell-j}} k_{s^{2m+2}}\right)_{s^{j-\left( m+1\right)}} k_{s^{\ell+m+1}} ds = \int P_4^{2m+2\ell+2}\left( k\right) ds \mbox{,}$$
where the highest derivative of $k$ that appears on the right is $k_{s^{m+\ell +1}}$.

Finally for each $j$ we look at the middle summation term of \eqref{E:intFin}.  Since $i \leq m < \ell$, the $k_{s^{m+i}}$ factors are no problem.  However, regardless of $i$, there will be a derivative factor of too high an order if $j>m$, so we need to do $j-m$ integrations by parts.  We have
$$\int k^2 k_{s^{m-i}} k_{s^{m+i}} k_{s^{\ell-j}} k_{s^{\ell+j}} ds = \left[ k^2 k_{s^{m-i}} k_{s^{m+i}} k_{s^{\ell-j}} k_{s^{\ell+j-1}} \right]_{\partial \gamma} - \int \left( k^2 k_{s^{m-i}} k_{s^{m+i}} k_{s^{\ell-j}} \right)_s k_{s^{\ell+j-1}} ds \mbox{.}$$
Clearly the boundary term here is equal to zero.  Integrating by parts a second time, if necessary, we have
\begin{multline*}
  \int k^2 k_{s^{m-i}} k_{s^{m+i}} k_{s^{\ell-j}} k_{s^{\ell+j}} ds\\
   = - \left[ \left( k^2 k_{s^{m-i}} k_{s^{m+i}} k_{s^{\ell-j}} \right)_s k_{s^{\ell+j-2}} \right]_{\partial \gamma} + \int \left( k^2 k_{s^{m-i}} k_{s^{m+i}} k_{s^{\ell-j}} \right)_{ss} k_{s^{\ell+j-2}} ds \mbox{.}
   \end{multline*}
If $\ell +j$ is odd then the boundary term is clearly equal to zero.  If $\ell + j$ is even then, as usual, we examine the other factor in the boundary term.  In this case it has the form $\left( e\, e\, o\, o\, e\right)_s$ or $\left( e\, e\, e\, e\, e\right)_s$.  More generally, with subsequent integrations by parts we are going to have $\left( e\, e\, o\, o\, e\right)_o$ or $\left( e\, e\, e\, e\, e\right)_o$.  Working similarly as in the proof of Lemma \ref{T:BD2} shows that such factors always have an $o$ factor, so these boundary terms are always equal to zero.  Another integration by parts, if necessary, gives
\begin{multline*}
  \int k^2 k_{s^{m-i}} k_{s^{m+i}} k_{s^{\ell-j}} k_{s^{\ell+j}} ds\\
   =  \left[ \left( k^2 k_{s^{m-i}} k_{s^{m+i}} k_{s^{\ell-j}} \right)_{ss} k_{s^{\ell+j-3}} \right]_{\partial \gamma} - \int \left( k^2 k_{s^{m-i}} k_{s^{m+i}} k_{s^{\ell-j}} \right)_{s^3} k_{s^{\ell+j-3}} ds \mbox{.}
  \end{multline*}
  Similar arguments as before give that the boundary term is again equal to zero.  
  
  In general, we have for each $i$,
  $$ \int k^2 k_{s^{m-i}} k_{s^{m+i}} k_{s^{\ell-j}} k_{s^{\ell+j}} ds = \left( -1\right)^{j-m} \int \left( k^2 k_{s^{m-i}} k_{s^{m+i}} k_{s^{\ell-j}} \right)_{s^{j-m}} k_{s^{\ell+m}} ds = \int P_6^{2m + 2\ell}\left( k\right) ds \mbox{,}$$
  with no higher derivative than $k_{s^{m+\ell}}$ appearing on the right hand side.\\
  
  Using all these results in \eqref{E:ksl2} we have for each $\ell>m$,
  $$\frac{d}{dt} \int k_{s^\ell}^2 ds = -2 \int k_{s^{m+\ell +2}}^2 ds + \int P_4^{2m+2\ell+2}\left( k\right) ds + \int P_6^{2m+2\ell}\left( k\right) ds \mbox{,}$$
  where no higher derivative that $k_{s^{m+\ell+1}}$ appears in the $\int P_4^{2m+2\ell+2}\left( k\right) ds$ term and no higher derivative than $k_{s^{m+\ell}}$ appears in $\int P_6^{2m+2\ell}\left( k\right) ds$.  On the second term we use \eqref{E:P4} while for the last we estimate using Proposition \ref{T:int} to estimate
  $$\int P_6^{2m + 2\ell}\left( k \right) ds \leq \varepsilon \int k_{s^{m+\ell+2}}^2 ds + c\left( \int k^2 ds\right)^{m+\ell+5} \mbox{,}$$
  where to have the second statement of Proposition \ref{T:int}, it is necessary to use that $P_6^{2m+2\ell}\left( k\right)$ contains derivatives of $k$ of order \emph{at most} $m+\ell +1$ (which is bigger than $m+\ell$, the highest order that actually occurs).  Therefore we obtain
  $$\frac{d}{dt} \int k_{s^\ell}^2 ds \leq - \left( \frac{\pi^2}{L_0^2}\right)^{m+2} \int k_{s^\ell}^2 ds + c\left( \int k^2 ds\right)^{m + \ell +5}$$
  from which the result follows since the second term above is bounded via Lemma \ref{T:PSW} and Lemma \ref{T:Length}. \hspace*{\fill}$\Box$
  
  \noindent \textbf{Remark:} As with the polyharmonic curve diffusion, we can establish the similar curvature blow up rate in the case that the maximal existence time $T$ is finite.  In this case, 
    $$\frac{d}{dt} \int k^2 ds = -2 \int k_{s^{m +2}}^2 ds + \int P_4^{2m+2}\left( k\right) ds + \int P_6^{2m}\left( k\right) ds \mbox{;}$$
  estimating the $P$ terms as above we get
  $$\frac{d}{dt} \int k^2 ds \leq c \left( \int k^2 ds\right)^{2m+5}$$
  from which the blow up rate again follows.\\[8pt]

\begin{kexp} \label{T:kexp}
If the initial curve $\gamma_0$ has sufficiently small scale-invariant energy
$$L^{2m+1} \int k_{s^m}^2 ds \leq c\left( m\right)$$
then, under \eqref{E:theflow} with normal speed \eqref{E:speed}, $\int k^2 ds$ decays exponentially.
\end{kexp}

\noindent \textbf{Proof:} We have using Lemma \ref{T:evlneqns}, the boundary conditions and \eqref{E:speed} 
\begin{equation} \label{E:k22}
  \frac{d}{dt} \int k^2 ds = 2\int \left( k_{ss} + \frac{1}{2} k^3 \right) \left[ \left( -1\right)^{m+1} k_{s^{2m+2}} + \sum_{j=1}^m \left( -1\right)^{j+1} k\, k_{s^{m+j}} k_{s^{m-j}} - \frac{1}{2} k\, k_{s^m}^2 \right] ds \mbox{.}
  \end{equation}
We will deal with each of these terms separately.  For the first term in \eqref{E:k22} we have by $m$ integrations by parts and using Lemma \ref{T:BCs}
$$\int k_{ss} k_{s^{2m+2}} ds = \left( -1\right)^m \int k_{s^{m+2}}^2 ds \mbox{.}$$
The next summation of terms in \eqref{E:k22} is
$$\sum_{j=1}^m\left( -1\right)^{j+1} \int k_{ss} k\, k_{s^{m+j}} k_{s^{m-j}} ds = \int P_4^{2m+2}\left( k \right) ds \mbox{.}$$
We want to verify that this can be rewritten such that the highest order derivative of $k$ is $k_{s^{m+1}}$.  For the $j$th term in the sum, this will require integrating by parts $j-1$ times.  For each integration by parts, the boundary term always contains a factor of an odd derivative of $k$, so the boundary terms are always equal to zero.  We obtain for each $j$
$$\int k_{ss} k\, k_{s^{m-j}} k_{s^{m+j}} ds = \left( -1\right)^{j-1} \int \left( k_{ss} k\, k_{s^{m-j}}\right)_{s^{j-1}} k_{s^{m+1}} ds = \int P_4^{2m+2}\left( k\right) ds \mbox{,}$$
where the highest order derivative of $k$ is $k_{s^{m+1}}$.

For the third term in \eqref{E:k22} we have
$$\int k_{ss} k\, k_{s^m}^2 ds = \int P_4^{2m+2}\left( k\right) ds$$
and the highest order derivative is $k_{s^m}$.

For the fourth term in \eqref{E:k22} we have similarly as the first by $m$ integrations by parts
$$\int k^3 k_{s^{2m+2}} ds = \left( -1\right)^m \int \left( k^3\right)_{s^{m+1}} k_{s^{m+1}} ds = \int P_4^{2m+2}\left( k\right) ds$$
with the highest order derivative of $k$ being $k_{s^{m+1}}$.

For the fifth summation term in \eqref{E:k22} we integrate the $j$th term by parts $j$ times observing again that every boundary term that arises contains an odd derivative of $k$ so is equal to zero.  We obtain
$$\int k^4 k_{s^{m-j}} k_{s^{m+j}} ds = \left( -1\right)^j \int \left( k^4 k_{s^{m-j}} \right)_{s^j} k_{s^m} ds = \int P_6^{2m}\left( k\right) ds$$
with highest order derivative of $k$ begin $k_{s^m}$.  The last term in \eqref{E:k22} already has the form $\int P_6^{2m}\left( k\right) ds$.

Therefore we have two kinds of term to estimate using Lemma \ref{T:int}.  We first use again \eqref{E:P4} together with Lemma \ref{T:PSW}:
$$\int P_4^{2m+2}\left( k \right) ds \leq c\, L \int k^2 ds \int k_{s^{m+2}}^2 ds \leq c\, L^{2m+1} \int k_{s^m}^2 ds \int k_{s^{m+2}}^2 ds \mbox{.}$$
Using the same interpolation inequalities we also have
$$\int P_6^{2m}\left( k \right) ds \leq c \left( \int k^2 ds\right)^2 \int k_{s^{m+1}}^2 ds \leq c\, L^{4m+2} \left( \int k_{s^m}^2 ds \right)^2 \int k_{s^{m+2}}^2 ds \mbox{.}$$

Substituting all these into \eqref{E:k22} we obtain
$$\frac{d}{dt} \int k^2 ds \leq \left[ -2 + c\, L^{2m+1} \int k_{s^m}^2 ds + c\, L^{4m+2} \left( \int k_{s^m}^2 ds \right)^2 \right] \int k_{s^{m+2}}^2 ds \mbox{.}$$
From Lemma \ref{T:Length} we know for sufficiently small initial energy that $L$ does not increase.  Also the energy itself is nonicreasing so if the above coefficient on the right hand side is initially less than $-\delta$ say, for some $\delta >0$, this will remain so and we have
$$\frac{d}{dt} \int k^2 ds \leq -\delta \int k_{s^{m+2}}^2 ds \leq -\tilde \delta\left( m, L_0\right) \int k^2 ds$$
hence the result.\hspace*{\fill}$\Box$

We may now prove similarly as for Corollary \ref{T:PHCDL3} exponential decay of curvature derivatives.

\begin{dexp} \label{T:dexp}
If the initial curve $\gamma_0$ has sufficiently small scale-invariant energy
$$L^{2m+1} \int k_{s^m}^2 ds \leq c\left( m\right)$$
then, under \eqref{E:theflow} with normal speed \eqref{E:speed}, $\int k_{s^{\ell}}^2 ds$ and $\left\| k_{s^\ell} \right\|_\infty$ decay exponentially for all $\ell$.
\end{dexp}

With these estimates in hand the exponential convergence of $\gamma$ to a unique horizontal straight segment now follows by exactly the same argument as in the previous section.  This completes the proof of Theorem \ref{T:main2}. \hspace*{\fill}$\Box$\\[8pt]

\end{document}